\documentclass{article}

\usepackage{amsmath,amssymb,amsthm,bm}
\usepackage{bbm}
\usepackage{enumerate}
\usepackage{setspace}
\usepackage[utf8]{inputenc}
\usepackage[dvipsnames]{xcolor}
\usepackage{tikz}
\usetikzlibrary[shadings]
\usepackage{mathrsfs}
\usepackage{url}
\usepackage{comment}

\makeatletter
\providecommand*{\diff}%
{\@ifnextchar^{\DIfF}{\DIfF^{}}}
\def\DIfF^#1{%
	\mathop{\mathrm{\mathstrut d}}%
	\nolimits^{#1}\gobblespace}
\def\gobblespace{%
	\futurelet\diffarg\opspace}
\def\opspace{%
	\let\DiffSpace\!
	\ifx\diffarg(%
	\let\DiffSpace\relax
	\else
	\ifx\diffarg[%
	\let\DiffSpace\relax
	\else
	\ifx\diffarg\{%
	\let\DiffSpace\relax
	\fi\fi\fi\DiffSpace}
\makeatother

\theoremstyle{plain}
\newtheorem{theorem}{Theorem}
\newtheorem*{theorem*}{Theorem}
\newtheorem{lemma}[theorem]{Lemma}
\newtheorem{proposition}[theorem]{Proposition}

\theoremstyle{definition}
\newtheorem{definition}[theorem]{Definition}

\theoremstyle{remark}
\newtheorem{remark}[theorem]{Remark}

\newcommand{\Z}{\mathbb{Z}}

\newcommand{\N}{\mathbb{N}}
\newcommand{\R}{\mathbb{R}}

\newcommand{\C}{\mathbb{C}}

\def\N{\mathbb{N}}
\def\Z{\mathbb{Z}}

\def\R{\mathbb{R}}
\def\C{\mathbb{C}}

\def\1{\mathbf{1}}

\newcommand{\norm}[1]{\|#1\|}

\title{A sharp discrete Sobolev inequality and left inverses of higher order differences}

\author{Raul Hindov, Shahaf Nitzan\footnote{The second author was partially supported by NSF CAREER grant DMS 1847796.},  Jan-Fredrik Olsen, Eskil Rydhe\footnote{The fourth author is supported by VR grant 2022-04307.}}

\begin{document}

\maketitle
\begin{abstract}
    For integers $\alpha,d \geqslant1$, we determine the best constants in the inequalities
    \[
    \norm{ (a_j)_{j=1}^{d} }_{\ell^1} \leqslant C_{\alpha,d} \norm{((\Delta^\alpha a)_j)_{j=1}^{d+\alpha}}_{\ell^2},
    \]
   where $\Delta$ is the difference operator. We find that these constants correspond to the $\ell^2$ norm of the discrete Chebyshev polynomials, and show that the extremal sequences may be obtained by appropriate sampling of the discrete Landau polynomials.  

From beginning to end, this project was carried out without use of AI.

\end{abstract}

\section{Introduction}

For finite sequences $a=(a_j)_{j=1}^d$ consider the following discrete inequalities of Poincar\`{e} type,
 \begin{equation}\label{general_discrete}
 \norm{ (a_j)_{j=1}^{d} }_{\ell^q} \leqslant C_{\alpha,q,p,d} \norm{((\Delta^\alpha a)_j)_{j=1}^{d+\alpha}}_{\ell^p},
 \end{equation}
where $p, q\in[1,\infty]$, $d,\alpha\in\N$. Here, $\Delta$ denotes the backward difference operator on sequences indexed by $\Z$, and we compute
	$(\Delta a)_j = a_j - a_{j-1}$
by extending  $(a_j)_{j=1}^d$ to a two sided-sequence, which is equal to zero on $\mathbb{Z} \backslash \{1,\ldots, d\}$.

The inequality \eqref{general_discrete} may be interpreted as a polynomial inequality: One may think of $(a_j)_{j=1}^d$ as the coefficients of a polynomial $p(z) = \sum_{j=1}^d a_j z^{j-1}$ of degree at most $d-1$,  and of $((\Delta^\alpha a)_j)_{j=1}^{d+\alpha}$ as the coefficients for the polynomial $(1-z)^{\alpha}p(z)$. Equivalently, putting $z=e^{it}$, one may interpret \eqref{general_discrete} as an inequality regarding the Fourier coefficients of the corresponding trigonometric polynomials.

A classical goal regarding such inequalities, is to find the optimal constant $C_{\alpha,q,p,d}$ in \eqref{general_discrete}, and determine the sequences for which equality is attained.

In the present paper, we provide a full solution for the case $q=1$, $p=2$, and all integer $\alpha \geqslant 1$.  To the best of our knowledge, this solution has not previously been determined. Our main result is as follows. 
\begin{theorem}\label{main_theorem}
	Let $\alpha, d \in \N$. For any $a_1,a_2,\dots,a_d \in \mathbb{C}$ the following inequality holds
	\begin{equation}\label{the_main}
		\norm{ (a_j)_{j=1}^{d} }_{\ell^1} \leqslant \sqrt{\binom{d+2\alpha}{2\alpha +1}\binom{2\alpha}{\alpha}^{-1}} \norm{((\Delta^\alpha a)_j)_{j=1}^{d+\alpha}}_{\ell^2}.
	\end{equation}
	This estimate is sharp, with equality obtained when $(a_j)_{j=1}^d$ are samples of the discrete Landau polynomial over the set $\{1,...,d\}$  (see \eqref{the L polynomial}).
\end{theorem}

We record some of the history in the study of the Poincar\`{e} type inequalities \eqref{general_discrete}. In the case $p=q=2$, the inequality \eqref{general_discrete} is also known as the discrete Wirtinger inequality. The best constant
 for $\alpha=1$ was obtained by Fan, Taussky and Todd, who also found the maximizing sequences \cite[Theorem 9]{FanTausskyTodd1955} (see also \cite{Izmestiev2017,Shisha1973}).
Higher-order discrete Wirtinger
inequalities have been studied by,   e.g., Milovanovi{\'c} and Milovanovi{\'c}
\cite{MilovanovicMilovanovic1997}.

Related extremal problems regarding trigonometric polynomials $t \mapsto p(e^{it})$ have also been studied, with the norm on the right hand side of \eqref{general_discrete} replaced by the $L^\infty(\mathbb{T})$ norm of $(1-z)^\alpha p(z)$, and the norm on the left hand side replaced by $|p(1)|$ (note that for non-negative coefficients, this quantity is equal to the $\ell^1$ norm of the Fourier coefficients). Under some additional restrictions, the corresponding constant was determined for $\alpha=1$ by Kravitz and Steinerberger \cite{Kravitz}, for $\alpha=2$ by Richardson \cite{Richardson2026}, and for $\alpha=3$ by Gaitán et al. \cite{Gaitan2026}.

Theorem \ref{main_theorem} is inspired by our work on similar problems for compactly supported functions on the real line. More specifically,  in \cite{Hindov}, we  prove the sharp inequality
\begin{equation} \label{continuous case}
 \|f\|_{L^1(-1,1)}
 \leq
 \frac{1}{(2\alpha-1)!!\sqrt{\alpha+\frac12}}\,
 \|f^{(\alpha)}\|_{L^2(-1,1)},
\end{equation}
for $f\in W_0^{\alpha,2}(-1,1)$, whose extremal functions are scalar multiples of the Landau kernels
$(1-x^2)^\alpha$.

The estimate \eqref{continuous case} was recently extended to all $p\in[1,\infty]$ by Garmanova and Kazimirov \cite{GarmanovaKazimirov2026}. For general $p, q \in [1,\infty]$, the analogous problem to \eqref{general_discrete} for compactly supported functions in $\R$ is far from solved in general. However, there are several elegant results; see, for example, the surveys in \cite{BoultonLang, Kuznetsov, Mitrinovic, NazarovShcheglova2025}.

Our proof of Theorem \ref{main_theorem} is similarly inspired by the argument we made when establishing \eqref{continuous case} in \cite{Hindov}. In particular, our approach revolves around a construction of explicit left-inverse operators of $\Delta^\alpha$, which allows us to reduce the problem of identifying the sharp constant to a norm-minimization problem over a class of discrete polynomials.  To establish Theorem \eqref{general_discrete}, below, we apply additional care with the location and multiplicity of the zeros of the corresponding polynomials. We first treat distinct integer zeros and then pass to arbitrary real zeros by  considering convex combinations.

\begin{remark}
The sharp continuous inequality \eqref{continuous case} can be recovered from Theorem~\ref{main_theorem} by standard discretisation techniques. As is well-known, under the corresponding rescaling and normalisation, the discrete Chebyshev polynomials of any fixed degree converge to the Legendre polynomials (see, e.g., \cite[p. 33]{Szego1975}). Similarly, the polynomials defining the extremal sequences in Theorem~\ref{main_theorem}, when suitably rescaled and normalised, converge to the Landau kernels $(1-x^2)^\alpha$.
\end{remark}

\section{Preliminary results}\label{section:prelim}

In this section we consider the difference operator $\Delta$ and  record some preliminary results that we will use in what follows.

\subsection{The difference operator $\Delta$}

The (backward) difference operator $\Delta$
is defined on complex valued sequences $(a_j)_{j \in \Z}$ by
\begin{equation} \label{discrete delta on Z}
	(\Delta a)_j = a_j - a_{j-1} 
    ,\qquad  j \in \Z.
\end{equation}

The following definition will allow us to formulate discrete analogues of results that are classical in the continuous case. While this material is not new, we provide short arguments to keep the text self-contained.

\begin{definition}
Let $(a_j)_{j\in\Z}$ be a real sequence. We say that $\mu \in \Z$ is a generalised zero for the sequence if one of the following holds:
\begin{enumerate}[(i)]
\item $a_{\mu}=0$,
\item $a_{\mu} \ \neq 0$ and $a_{\mu} a_{\mu+1}<0$.
\end{enumerate}
\end{definition}

The following lemma can be considered a discrete analogue of Rolle's theorem.
\begin{lemma}[Discrete Rolle's lemma]\label{lemma:DiscreteRolle}
    Suppose $(a_j)_{j \in \Z}$ is a sequence of real numbers that has generalised zeroes at indices $\mu, \nu \in \Z$ with $\mu < \nu$. Then the sequence $(\Delta a_j)_{j\in \Z}$ has a generalised zero for some $j \in \{\mu+1, \ldots, \nu\}$.
\end{lemma}
\begin{proof}

Assume for a contradiction that $(\Delta a_j)_{j={\mu+1}}^{\nu}$ has no generalized zeroes. Without loss of generality, we may assume that $\Delta a_j > 0$ for all $j \in \{\mu +1, \ldots, \nu\}$ and that $\Delta a_{\nu+1} \geq 0$.  Since $a_{\mu}$ is a generalized zero, and $\Delta a_{\mu+1}>0$, this implies that $a_{\mu+1}> 0$. It follows that the sequence $(a_j)_{j={\mu+1}}^{\nu}$ is an increasing sequence of positive numbers. Since $\Delta a_{\nu+1}\geq 0$, we have $0 < a_{\nu}\leq a_{\nu+1}$, which contradicts the assumption that $a_{\nu}$ is a generalized zero. 
\end{proof}

\subsection{The polynomials $r_{\alpha}$ and functions $f_{\alpha}$}

We define a family of polynomials $\{r_\alpha\}_{\alpha \in \N}$ by
\begin{equation}\label{r_a}
	r_\alpha(x):=\frac{1}{\alpha !}\prod_{j=1}^\alpha(x+j),\qquad \alpha\in\N,x\in\R.
\end{equation}
It will be convenient notationally to let $r_0(x):=1$.

We record the following useful fact as a lemma.
\begin{lemma} \label{r recursion lemma}
    For all $\alpha \in \N$, we have $r_\alpha(0) = 1$ and
    \begin{equation} \label{Delta r_alpha}
      \Delta r_\alpha = r_{\alpha-1}.
    \end{equation}
\end{lemma}
\begin{proof}
    This follows by direct application of the definition of $r_\alpha$.
\end{proof}

Since $r_{\alpha}(0)=1=r_{\alpha+1}(0)$, the identity (\ref{Delta r_alpha}) allows us to obtain the  summation formula
\begin{equation}\label{sum r_a}
\begin{aligned}
\sum_{j=0}^{k}r_{\alpha}(j) =r_{\alpha+1}(0)+\sum_{j=1}^{k}\Delta r_{\alpha+1}(j)
=r_{\alpha+1}(k).
\end{aligned}
\end{equation}
Next, consider the functions $f_{\alpha}:\Z\rightarrow\Z$ defined by
  \begin{equation}\label{f_a}
        f_\alpha(n) = \left\{ \begin{aligned}
        r_{\alpha-1}(n) & \quad n \geq  -\alpha+1 \\
        0 & \quad n \leq -1.
        \end{aligned}\right.
    \end{equation}
    Observe that $f_\alpha$ is well defined since $r_{\alpha-1}(n) = 0$ for $n \in \{-\alpha+1,\ldots, -1\}$. Moreover, as $r_\alpha(0) = r_{\alpha-1}(0)=1$, Lemma \ref{r recursion lemma} implies that
    \begin{equation}\label{Df_a}
    \Delta f_\alpha = f_{\alpha-1}.
    \end{equation}

	Over the next two lemmas, we study the (generalized) zero and sign structure of functions of the form $f_{\alpha}-p$, where $p$ is a polynomial of degree at most $\alpha-1$.
	
	\begin{remark}
	We use the subscript $\alpha$ for $f_\alpha$, even though it is constructed from two polynomials of degree $\leq \alpha-1$, since, as can be seen in the next three lemmas, the zero structure of $f_\alpha-p$ mimics that of a polynomial of degree $\alpha$.
	\end{remark}

\begin{lemma}[Discrete broken line lemma] \label{Discrete broken line lemma}
	For $\alpha\in \N$, let $p$ be a polynomial of degree at most $\alpha-1$.  If $g(n) = f_\alpha(n) - p(n)$ has at least $\alpha+1$ generalised zeroes, then either $p(n) = r_{\alpha-1}(n)$ for all $n \in \Z$ or $p(n) = 0$ for all $n \in \Z$.
\end{lemma}

\begin{proof}
The proof is by induction. First, let $\alpha=1$. In this case, $p$ is a constant function and
\begin{equation*}
    f_1(n) = \left\{
    \begin{aligned}
        1 & \quad n \geq 0 \\
        0 & \quad n \leq -1
    \end{aligned}\right..
\end{equation*}
In particular, for $g = f_1 - p$  to have at least $2$ distinct generalised zeroes, we must have either $ p(n) = 1$ for all $n \in \Z$ or $ p(n) = 0$ for all $n \in \Z$.

Next, let $\alpha \geq  2$ and suppose that $g = f_\alpha - p$ has at least  $\alpha+1$ distinct generalised zeroes. By Lemma \ref{lemma:DiscreteRolle}, we obtain that $\Delta g = \Delta f_{\alpha} - \Delta p$ has at least $\alpha$ distinct generalised zeroes. Since $\Delta p$ is a polynomial of degree at most $\alpha-2$, we may apply (\ref{Df_a}) and the induction assumption to conclude that $\Delta p$ is either identically equal to $r_{\alpha-2}$ or identically equal to $0$ on $\Z$.  

By \eqref{Delta r_alpha}, this implies that there exists a constant $C \in \R$ so that one of the following holds:
\begin{align*}
    (i)& \quad p(n) = r_{\alpha-1}(n) + C \quad \forall n \in \Z, \\[2mm]
    (ii)& \quad p(n) = C \quad \forall n \in \Z.
\end{align*}

If (ii) holds, it suffices to observe that, for $\alpha >  1$, the function $f_\alpha$ is strictly monotone  over  $n>0$. This implies that if $C \neq 0$, then $g = f_\alpha - p$ has at most one generalised zero. Hence, we must have $C=0$. If (i) holds with $C\neq 0$, then it follows immediately that $f_\alpha$ has no generalised zeroes for $n \geq -\alpha+1$, and, since $r_{\alpha-1}$ is strictly monotone for $n \leq -\alpha+1$, it has at most one generalised zero for $n \leq	 -\alpha$. Therefore $C=0$, which ends the proof.
\end{proof}

\begin{lemma}\label{lemma sign structure}
  For $\alpha\in \N$,  let $p$ be a polynomial of degree at most $\alpha-1$.  If $g(n) = f_\alpha(n) - p(n)$ has at least $\alpha$ generalised zeroes $k_1 < k_2 < \ldots < k_{\alpha}$, then
  \begin{equation*}
  	(-1)^{\alpha+\ell} g(n) \geq 0 \qquad \forall  n  \in \{k_\ell+1,..., k_{\ell+1}\}
  \end{equation*}
where $\ell\in \{0,1,\ldots, \alpha\}$ with 
 $\{k_0+1, \ldots, k_{1}\} = \{n : n \leq k_1\}$ and $\{k_\alpha+1, \ldots, k_{\alpha+1}\} = \{n: n > k_\alpha\}$.
\end{lemma}
\begin{proof}
If $g$ has at least $\alpha+1$ generalized zeroes, the conclusion follows from Lemma \ref{Discrete broken line lemma} since either $g = f_\alpha - 0 = f_\alpha$ or $$g = f_\alpha - r_{\alpha-1} =  \left\{ \begin{aligned} 0 & \quad \text{if }\ n \geq - \alpha+1 \\ -r_{\alpha-1} & \quad \text{if }\ n \leq -1 \end{aligned} \right.$$ In either case, the conclusion follows trivially. 

We now turn to the case when $g$ has exactly $\alpha$ zeroes $k_1 < \cdots < k_\alpha$.
The proof is by induction. First, let $\alpha=1$. In this case $p$ is a constant function and
\begin{equation*}
    f_1(n) = \left\{
    \begin{aligned}
        1 & \quad n \geq 0 \\
        0 & \quad n \leq -1
    \end{aligned}\right.
\end{equation*}
So, if $g = f_1 - p$ has \textit{exactly} one generalised zero, $k_1$, then we must have $p(n) = c$ with $0<c<1$ and $k_1=-1$. This implies the statement of the lemma.

Next, assume $\alpha \geq 2$ and suppose that $g = f_\alpha - p$ has exactly $\alpha$ distinct generalised zeroes. By Lemma \ref{lemma:DiscreteRolle}, we obtain that $\Delta g =   f_{\alpha-1} - \Delta p$ has at least $\alpha-1$ distinct generalised zeroes.

By Lemma \ref{Discrete broken line lemma}, it follows that $\Delta g$ cannot have more than $\alpha-1$ distinct generalised zeroes. Indeed, if this were to happen, then we would have that $\Delta p=0$ or  $\Delta p=r_{\alpha-2}$, which would imply that $g$ has either an infinite amount of zeroes or at most one generalized zero (cf., the last paragraph of the proof of Lemma \ref{Discrete broken line lemma}, above). In both cases this contradicts our assumption.

We may therefore conclude that $\Delta g$ has exactly $\alpha-1$ distinct generalised zeroes $\nu_1 < \cdots < \nu_{\alpha-1}$. Moreover, Lemma \ref{lemma:DiscreteRolle} implies that they interlace with the generalized zeroes of $g$, that is,
\begin{equation} \label{so many labels}
	k_\ell < \nu_\ell  \leq k_{\ell+1}\qquad \ell \in \{ 1, \ldots, \alpha-1\}.
\end{equation}

Fix  $\ell \in \{1, \ldots, \alpha\}$. By the induction hypothesis (and rewriting the interval), it follows that
  \begin{equation} \label{label label}
  	(-1)^{\alpha+\ell} \Delta g(n) \geq 0 \qquad \forall n  \in \{\nu_{\ell-1}+1,..., \nu_{\ell}\} ,
  \end{equation}
  where $\{\nu_0+1, \ldots, \nu_1\}=\{n : n \leq \nu_1\}$ and $\{\nu_{\alpha-1}+1, \ldots, \nu_{\alpha}\}=\{n : n > \nu_{\alpha-1}\}$.
Since $k_\alpha\geq \nu_{\alpha-1}$, we have that  (\ref{label label}) holds for $n=k_\alpha +1$. For $\ell \in \{1, \ldots, \alpha-1\}$, one can check that (\ref{so many labels}) implies (\ref{label label}) for $n=k_{\ell}+1$. Explicitly this means that
  \[
  (-1)^{\alpha+\ell} g(k_{\ell}+1)\geq (-1)^{\alpha+\ell} g(k_{\ell}) \qquad \forall \ell \in \{1, \ldots, \alpha\}.
  \]
  As each $k_{\ell}$ is a generalized zero for $g$ we may conclude that
   \[
  (-1)^{\alpha+\ell} g(k_{\ell}+1)\geq 0.
  \]
  Since $g$ has no sign changes over the interval $\{k_{\ell}+1,...,k_{\ell+1}\}$ the same conclusion holds for every $k$ in this interval. This completes the proof for  $\ell \in \{1, \ldots, \alpha\}$.

The case $\ell=0$ may be treated in a similar way: We note that since $k_1<\nu_1$ the function $g$ is monotone over  $\{n : n \leq k_1\}$ and its only generalized zero over this interval is at the point $k_1$. This implies that if $g$ decreases over this interval then it must be positive there, whereas if it increases then it must be negative. Applying the induction hypothesis we obtain the statement for the case  $\ell=0$.
\end{proof}

We conclude this section with the following observation.
\begin{lemma}\label{sum f_a}
Fix  $\alpha\in\N$. For $j\in\Z$ let $p_j$ be a polynomial of degree at most $\alpha-1$ and denote
\[
g_j(k)=f_{\alpha}(j-k)+p_j(k) \qquad k\in\Z.
\]
Then, for every $m\in\N$, the function $\psi_m$ defined by
\[
\psi_m(k) := \sum_{j=-m}^{m} g_j(k)
\]
is equal over the interval $\{-m,...,m\}$ to a polynomial of degree $\alpha$ and with leading coefficient $(-1)^{\alpha}/(\alpha!)$.
\end{lemma}
\begin{proof}
Indeed, to see this first note that
\[
\psi_m(k) = \sum_{j=-m}^{m} \Big( f_\alpha(j-k)+ p_j(k) \Big) =\sum_{j=-m}^{m} f_\alpha(j-k)+P(k),
\]
where $P$ is a polynomial of degree at most $\alpha-1$. Next, observe that for $k\in\{-m,...,m\}$ we have
\[
\sum_{j=-m}^{m} f_\alpha(j-k)=\sum_{j=k}^{m} r_{\alpha-1}(j-k)=\sum_{j=0}^{m-k} r_{\alpha-1}(j).
\]
Applying the identity (\ref{sum r_a}), we find that over the interval $\{-m,...,m\}$ we have
\[
\psi_m(k) = r_{\alpha}(m-k)+P(k).
\]
The conclusion follows.
\end{proof}

\section{Proof of Theorem \ref{main_theorem} (part 1 of 2)}\label{section:Proof}

In this section we prove the estimate
\begin{equation}\label{poly estimate}
\norm{ (a_j)_{j=1}^{d} }_{\ell^1} \leqslant
\frac{1}{\alpha!}
 \min_{p\in\mathscr{P}_\alpha} \left(\sum_{k=1}^{d+\alpha}|p(k)|^2\right)^\frac12
\norm{((\Delta^\alpha a)_j)_{j=1}^{d+\alpha}}_{\ell^2}
\end{equation}
where $\mathscr{P}_\alpha$ is the family of polynomials with leading coefficient equal to 1 and of degree equal to $\alpha$.

\subsection{Preliminaries on difference operators $\Delta^\alpha$ and $\Delta^\alpha_d$}

We first fix some notations.  Given $d,\ell\in\N$ we let  $\mathbb{C}^d$ be the space of complex valued $d$-tuples and $\mathbb{C}^{d\times \ell}$ denote the space of complex matrices with $d$ rows and $\ell$ columns. To distinguish, we denote the $j$'th entry of a $d$-tuple $a\in \mathbb{C}^d$ by $a_j$ or $(a)_j$, whereas for a matrix $A\in \mathbb{C}^{d\times \ell}$ we use the notation $A(j,k)$ for the entry of $A$ at row $j$ and column $k$.
We identify $A\in\mathbb{C}^{d\times \ell}$ with the linear operator $\mathbb{C}^{\ell}\to\mathbb{C}^d$, $x\mapsto Ax$. We denote the adjoint of $A$, which for real valued matrices coincides with the transpose of $A$, by $A^*$.

Fix $d, \alpha \in\N$. We embed $(a_j)_{j=1}^d \in \mathbb{C}^d$ into the space of two-sided sequences (e.g., indexed by $\Z$) by setting $a_j = 0$ for all $j \in \Z\setminus  \{1,2, \ldots, d\}$. Applying this embedding, we find that the difference operator  $\Delta$, defined in \eqref{discrete delta on Z}, induces a transformation $\Delta_d:\mathbb{C}^d\rightarrow \mathbb{C}^{d+1}$ given by 
\begin{equation*}
    (\Delta_d a)_j = \left\{ \begin{aligned}
    a_1 \quad & \qquad j=1 \\
    a_j - a_{j-1} & \qquad 2 \leq j \leq d \\
    - \ a_d \quad & \qquad j=d+1.
    \end{aligned}\right.
\end{equation*}
We let $\Delta_d^\alpha: \mathbb{C}^d \rightarrow \mathbb{C}^{d+\alpha}$ denote the operator defined by
\[
\Delta_d^\alpha:= \Delta_{d+\alpha-1} \ldots \Delta_{d+1} \Delta_d.
\]
We remark that $\Delta_d^\alpha$ may be obtained from $\Delta^\alpha$ via the natural embedding described above. 
The adjoint of $\Delta_d$ is the operator $\Delta_d^\ast: \mathbb{C}^{d+1}\rightarrow \mathbb{C}^d$ defined by
\begin{equation*}
    (\Delta_d^\ast b)_j = b_{j}-b_{j+1}, \qquad j \in \{1,2,\ldots, d\}.
\end{equation*}
A basic observation, which we will repeatedly use, is that matrix multiplication by $\Delta_d$ from the right has a useful interpretation in terms of the relation $B\Delta_d = (\Delta_d^\ast B^\ast)^\ast$. Explicitly, for a matrix $B: \mathbb{C}^{d+1} \rightarrow \mathbb{C}^d$, with  $j \in \{1,2,\ldots, d\}$, $k \in \{1,2, \ldots, d \}$, we have
\begin{align*}
	(B\Delta_d)(j,k) &= \big( \Delta_d^\ast B(j,\cdot) \big)(k) \\[2mm] &= B(j,k) - B(j,k+1).
\end{align*}
This observation can be extended to $\Delta_d^\alpha$, with $\alpha \geq 2$, and matrices  $B: \mathbb{C}^{d+\alpha} \rightarrow \mathbb{C}^d$,  by
\begin{equation} \label{explicitly}
    (B\Delta_d^\alpha)(j,k) =\big( (\Delta_d^\alpha)^* B(j,\cdot) \big)(k). 
\end{equation}

The next lemma characterises the kernel of the operator $(\Delta^\alpha_d)^\ast$.

\begin{lemma}\label{prop:KernelOfDeltaStar}
	Let $\alpha,d\in\N$, $\textcolor{black}{b}\in \mathbb{C}^{d+\alpha}$. The following are equivalent:
	\begin{enumerate}[$(i)$]
		\item $(\Delta^{\alpha}_d)^*  \textcolor{black}{b} = 0$.
		\item There exists a polynomial $p$ of degree $\leqslant\alpha-1$ such that $b_k=p(k)$ for $k\in\{1,\ldots,d+\alpha\}$.
	\end{enumerate}
\end{lemma}

\begin{proof}
Assume $(ii)$ and recall that $(\Delta_d^\alpha)^\ast = \Delta^\ast_{d}   \ldots   \Delta^\ast_{d+\alpha-1}$.
Then
\[
(\Delta^*_{d+\alpha-1} b)(k)= b_k-b_{k+1}
\]
 depends polynomially on $k\in\{1,\ldots,d+\alpha-1\}$ with degree at most $\alpha-2$. The statement   (i) now follows by induction. 
	
	\textcolor{black}{Next, assume (i).} The monomial sequences $(k^j)_{k=1}^{d+\alpha}$, $0\leqslant j\leqslant \alpha-1$, are linearly independent as vectors in $\mathbb{C}^{d+\alpha}$. In light of the first part of the proof, it follows that they generate an $\alpha$-dimensional subspace of $\ker (\Delta^\alpha_d)^\ast$. But $(\Delta^\alpha_d)^\ast \colon \mathbb{C}^{d+\alpha}\to\mathbb{C}^d$ has full rank because $\Delta^\alpha_d\colon \mathbb{C}^d\to\mathbb{C}^{d+\alpha}$ is injective. That is, $\dim\ker (\Delta^\alpha_d)^* = \alpha$. We conclude that the sequences $(k^j)_{k=1}^{d+\alpha}$, $0\leqslant j\leqslant \alpha-1$, form a basis for $\ker(\Delta^\alpha_d)^*$, so $(i)$ implies $(ii)$.
\end{proof}

\subsection{On the structure of the left inverses of $\Delta_d^\alpha$}

Let $\alpha,d\in\N$ be fixed. We are interested in the structure of left-inverses of
$\Delta_d^\alpha$, i.e., matrices $B_{\alpha}\in\mathbb{C}^{d + \alpha} \rightarrow \mathbb{C}^d$ such that
\begin{align}\label{eq:EquationForLeftInverse}
	B_{\alpha}\Delta_d^\alpha = I.
\end{align}
Note that we suppress the dependence of $B_\alpha$ on the parameter $d$.

Over the next few claims, we study the general solution of the linear equation \eqref{eq:EquationForLeftInverse}. First, in Proposition  \ref{prop:hiho}, we consider the corresponding homogeneous equation. Next, we find a particular solution to  \eqref{eq:EquationForLeftInverse}, and finally in Lemma \ref{lemma:StructureOfLeftInverses} we describe the general solution to this equation.

\begin{proposition}\label{prop:hiho}
	For a matrix $B_{\alpha}^{(h)}: \C^{d+\alpha} \rightarrow \mathbb{C}^d$, the following are equivalent:
	\begin{enumerate}[$(i)$]
		\item $B_{\alpha}^{(h)}\Delta^\alpha_d =0$.
		\item For each $j=1,...,d$ there exists a polynomial $p_j$ with degree $\leqslant \alpha-1$ such that $B_{\alpha}^{(h)}(j,k)=p_j(k)$.
	\end{enumerate}
\end{proposition}
\begin{proof} Keeping in mind \eqref{explicitly}, this follows immediately from Lemma~\ref{prop:KernelOfDeltaStar}.
\end{proof}

To express a particular solution to \eqref{eq:EquationForLeftInverse}, recall the family $\{f_\alpha\}_{\alpha \in \N}$ defined in (\ref{f_a}) and consider the matrix
\begin{equation}\label{B_a}
	B_{\alpha}^{(p)}(j,k)=f_{\alpha}(j-k),
\end{equation}
where $j \in \{1, \ldots, d\}$ and $k \in \{1, \ldots, d+\alpha\}$.

We have the following result on left inverses of the matrices $\Delta_d^\alpha$.

\begin{proposition}\label{prop:silver}
	Let $\alpha, d \in \N$, and $B_{\alpha}^{(p)}: \mathbb{C}^{d+\alpha} \rightarrow \mathbb{C}^d$ be as in (\ref{B_a}).
	Then, for $\Delta_d^\alpha= \Delta_{d+\alpha-1} \ldots \Delta_{d+1} \Delta_d$, the following hold:
    \begin{enumerate}[(i)]
        \item $B_{\alpha}^{(p)}\Delta_{d+\alpha-1} = B_{\alpha-1}^{(p)}$.
        \item $B_{\alpha}^{(p)}\Delta_d^\alpha = I_d$.
    \end{enumerate}
\end{proposition}
\begin{proof}
	As in the proof of Proposition \ref{prop:hiho},
	  we use \eqref{explicitly} to obtain
    \begin{equation*}
        (B_\alpha^{(p)}\Delta_{d+\alpha-1})(j,k) = \big( (\Delta_{d+\alpha-1})^\ast B_\alpha^{(p)}(j,\cdot) \big)(k).
    \end{equation*}
     Since $B_\alpha^{(p)}(j,k) = f_\alpha(j-k)$, the identity (\ref{Df_a}) implies that
     \begin{align*}
        (\Delta_{d+\alpha-1})^\ast B_{\alpha}^{(p)}(j,k)
        &=
         f_{\alpha}(j-k) - f_{\alpha}(j-k-1) \\[2mm]
       &=
        f_{\alpha-1}(j-k)
        = B_{\alpha-1}^{(p)}(j,k).
     \end{align*}
    This establishes (i). Statement (ii) follows  by induction.
\end{proof}

An immediate consequence of propositions \ref{prop:hiho} and \ref{prop:silver} is the following structure lemma for left inverses of $\Delta_d^\alpha$.

\begin{lemma}\label{lemma:StructureOfLeftInverses}
	For a matrix $B_\alpha: \C^{d+\alpha} \rightarrow \mathbb{C}^d$, the  following are equivalent.
	\begin{enumerate}[$(i)$]
		\item $B_\alpha$ is a left inverse of $\Delta_d^\alpha$.
		\item For each $j \in \{1,\ldots, d\},$ there exists a polynomial $p_j$, of degree at most  $\alpha-1$, such that for $k \in \{1, \ldots, d+\alpha\}$ we have $$B_\alpha(j,k)=B_\alpha^{(p)}(j,k)-p_j(k).$$
	\end{enumerate}
\end{lemma}

The freedom which Lemma \ref{lemma:StructureOfLeftInverses} gives in choosing the row polynomials will allow us to control the sign and size of the sums of the columns of $B_\alpha$, as we do this in the next subsection.

\subsection{Left inverses with prescribed column signs and sums}

The following definitions introduce terminology that we use to formulate Theorem \ref{main thm}, below.
\begin{definition}
For $S\subseteq \Z$, consider the sequences $a=(a_n)_{n\in S}$ and $b=(b_n)_{n\in S}$. We say that $b$ \textit{follows the sign structure} of $a$ (over $S$) if the following two conditions hold.
\vspace{5pt}
\begin{itemize}
\item[(i)] For every $n\in S$ we have $a_nb_n\geq 0$.
\vspace{2pt}
\item[(ii)] If for some $n\in S$ we have $a_n=0$ then we also have $b_n=0$.
\end{itemize}
\end{definition}

Note that while \textit{following the sign structure} is not a symmetric relation, it is a transitive one. That is, if $c$ follows the sign structure of $b$, and $b$ follows the sign structure of $a$, then $c$ follows the sign structure of $a$.

\begin{definition}\label{def poly realization}
Let $p$ be a real valued polynomial. We say that a matrix $A: \mathbb{C}^m \rightarrow \mathbb{C}^d$  \textit{realizes the polynomial} $p$ if
\vspace{5pt}
\begin{itemize}
\item[(i)] For every $j \in \{1, \ldots, d\}$ we have that   $\{A(j,k)\}_{k=1}^{m}$ follows the sign structure of $\{p(k)\}_{k=1}^{m}$ over $S=\{1,\ldots,m\}$.
\vspace{5pt}
\item[(ii)]  For every $1\leq k\leq m$ we have
\[
\sum_{j=1}^d A(j,k)=p(k).
\]
\end{itemize}
\end{definition}

\begin{remark}\label{abs value equality}
Observe that if a matrix $A$ realizes a polynomial $p$ in this way, then  for $k \in \{1, \ldots, m\}$, we have
\[
\sum_{j=1}^d |A(j,k)|=\Big|\sum_{j=1}^d A(j,k)\Big|=|p(k)|.
\]
This property will be used in the next subsection.
\end{remark}

We are now ready to prove the following theorem.
\begin{theorem}\label{main thm}
Fix $d,\alpha\in\N$, and let $p$ be a real polynomial of degree $\alpha$. Assume that the leading coefficient of $p$ is $(-1)^{\alpha}/\alpha!$ and that the zeroes of $p$ are all real. Then there exists a left inverse of $\Delta_d^\alpha$ which realizes $p$.
\end{theorem}
\begin{proof}
We divide the proof into three steps. First we show that the claim holds in the case where the zeroes of $p$ are distinct integers, next we remove the requirement that these integers are distinct, and finally, in Step III, we consider the case of general real zeroes.

\vspace{5pt}

\underline{\textit{Step I:}}
Let
\begin{equation*}
	p(x)= \frac{1}{\alpha!} \prod_{n=1}^{\alpha}(k_n-x),
\end{equation*}
	 where $k_n$ are distinct integers. Denote by $\{\mathcal{L}_n\}_{n=1}^{\alpha}$   the Lagrange interpolation basis with respect to these points, that is,
	\begin{align*}
		\mathcal{L}_n(x)=\prod_{\substack{1\le \ell\le \alpha;\\ \ell\ne n}}\frac{k_\ell-x}{k_\ell-k_n}.
	\end{align*}
	For $j, k \in \Z,$ define
	\begin{equation}
	\begin{aligned} \label{the definition}
		g_j(k) &=f_\alpha(j-k) - \sum_{n=1}^\alpha \mathcal{L}_n(k)f_\alpha(j-k_n) \\[2mm]
		&:= f_\alpha(j-k)+s_j(k),
	\end{aligned}
	\end{equation}
and observe that  $s_j(k)$ are polynomials of degree at most $\alpha-1$ in $k$.

Next, we claim that the functions $g_j$ satisfy the following properties:

\vspace{2pt}

\begin{itemize}
\item[(i)] For every $j\in\Z$ the set $\{k_n\}_{n=1}^{\alpha}$ belongs to the zero set of $g_j$. In particular, this implies that each of the functions $g_j$ has at least $\alpha$ zeroes.

\vspace{0pt}

\item[(ii)] For every fixed $j\in\Z$, the sequence $\{g_j(k)\}_{k\in\Z}$ follows the sign structure of $\{p(k)\}_{k\in\Z}$.

\vspace{0pt}

\item[(iii)] For every $j\leq 0$ and for every $j\geq d+1$, the function $g_j$ is equal to a polynomial of degree at most $\alpha -1$ over the interval $\{1,...,d+\alpha\}$.

\vspace{0pt}

\item[(iv)]  For every $m\in \N$ the sum
\[
\sum_{j=-m}^{m} g_j(k)
\]
is equal over $\{-m,...,m\}$ to a polynomial of degree $\alpha$ and with leading coefficient $(-1)^{\alpha}/(\alpha !)$.

\end{itemize}

Indeed, (i) follows from the definition of the Lagrange interpolation basis. Next, (ii) follows by keeping in mind property (i) and applying Lemma \ref{lemma sign structure}. For (iii), we recall the definition of $f_{\alpha}$: If $j\leq 0$ then $f_{\alpha}(j-k)=0$ over $\{1,...,d+\alpha\}$, while if $j\geq d+1$ then $f_\alpha(j-k)=r_{\alpha-1}(j-k)$ over this interval. Finally, (iv) follows from Lemma \ref{sum f_a}.

We are now ready to construct the required left inverse. Put $\rho=\max \{|k_n|: n=1,..., \alpha\}$ and fix $m\in\N$ so that $m\geq \max\{d+\alpha, \rho\}$. Denote
\[
G(k)=\sum_{j=-m}^0g_j(k)+\sum_{j=d+1}^m g_j(k),\qquad k\in\Z,
\]
and note that properties (i) and (ii) hold for $G$ in place of $g_j$. Moreover, observe that, by Property (iii), it follows that $G$ is a polynomial of degree at most $\alpha-1$ over the interval $\{1,...,d+\alpha\}$.

We define the first row of the left inverse by declaring
\begin{align*}
B(1,k)&=g_1(k)+G(k) \\[2mm] &=f_\alpha(1-k)+s_1(k)+G(k) \qquad k\in\{1,...,d+\alpha\}.
\end{align*}
We next define all other rows of the matrix by putting, for $j=2,...,d$,
\begin{align*}
B(j,k) &=g_j(k) \\[2mm] &=f_\alpha(j-k)+s_j(k) \qquad k\in\{1,...,d+\alpha\}.
\end{align*}
Keeping in mind that the functions $s_j$ are polynomials of degree at most $\alpha-1$, and that over the interval $\{1,...,d+\alpha\}$ the same is true for $G$, we may apply Lemma \ref{lemma:StructureOfLeftInverses} to conclude that $B(j,k)$ is indeed a left inverse of $\Delta_d^{\alpha}$.

Moreover, due to property (ii), and since this property holds also for $G$, we find that for every $j=1,...,d$ the sign structure of $\{B(j,k)\}_{k=1}^{d+\alpha}$ follows the sign structure of  $\{p(k)\}_{k=1}^{d+\alpha}$.

It remains to show that the sum of the rows is as required. For this, note that due to the definition of $G$ we have
\[
\sum_{j=1}^{d}B(j,k)=\sum_{j=-m}^{m}g_j(k)\qquad k \in\{1,...,d+\alpha\}.
\]
The sum on the right hand side may be considered as a function over $\Z$. Property (iv) implies that over the interval $\{-m,...,m\}$ this function is equal to some polynomial $q$ of degree $\alpha$ and
with leading coefficient $(-1)^{\alpha}/(\alpha !)$. The zero set of the polynomial $q$ is equal to $\{k_n\}_{n=1}^\alpha$. Indeed, since $m\geq |k_n|$ for every $n\in\{1,...,\alpha\}$, and $q$ is a polynomial of degree $\alpha$, this is implied by property (i), above. 

Finally, noting that the polynomials $p$ and $q$ are of the same degree, and that they have the same zero set and the same leading coefficient, we may conclude that $q=p$. This completes the proof for the first step.

\vspace{5pt}

\underline{\textit{Step II:}} Assume that the zeroes of $p$ are  integers, but not necessarily distinct. This means that $p$ has a decomposition
\[
p(x)=\frac{(-1)^{\alpha}}{\alpha!}\prod_{n=1}^{K}(x-k_n)^{\beta_{n}},
\]
with $k_1, \ldots, k_K\in\Z$ distinct. We prove by induction on the \textit{multiplicity degree} of $p$, which we define by
\[
\mathrm{md}(p):=\alpha-K=\sum_{n=1}^{K}\beta_{n}-K.
\]
The case $\mathrm{md}(p)=0$ was proved in Step I above. So, assume that $\mathrm{md}(p)\geq 1$. This means that there exists $n_0 \in \{1, \ldots, K\}$ so that $\beta_{n_0}\geq 2$.
For this fixed value of $n$, put
\[
v :=\min\{t\in\Z: k_{n_0}<t,\: p(t)\neq 0\}
\]
and
\[
\quad
u :=\max\{t\in\Z: t<k_{n_0},\:p(t)\neq 0\}.
\]
Denote
\[
p_u(x):=\frac{x- u}{x-k_{n_0}}p(x) \quad \textrm{ and } \quad \quad p_v(x):=\frac{x- v}{x-k_{n_0}}p(x).
\]
Observe that the multiplicity degrees of $p_u$ and $p_v$ satisfy $$\mathrm{md}(p_u)=\mathrm{md}(p_v)=\mathrm{md}(p)-1.$$ By the induction assumption, there exist left inverses of $\Delta_d^\alpha$ which realize $p_{u}$ and $p_{v}$. Denote them by $B_{\alpha}^{u}$ and $B_{\alpha}^{v}$ respectively.

It is straight-forward to check that both $p_{u}$ and $p_{v}$ follow the sign structure of $p$ over $\Z$. We may therefore apply the transitivity property to conclude that for every $j=1,..., d$, the sign structures of $B_{\alpha}^u(j,\cdot)$ and $B_{\alpha}^v(j,\cdot)$ follow the sign structure of $p$ over  $\{1,\ldots,d+\alpha\}$.

Next, as $u<k_{n_0}<v$  there exists $0<\tau<1$ so that $k_{n_0} = \tau u  +(1-\tau)v$. Put
\[
B_{\alpha}:=\tau B_{\alpha}^{u}+(1-\tau)B_{\alpha}^{v}.
\]
As $\tau$ is positive, we find that for every $1\leq j\leq d$, the sign structure of $B_{\alpha}(j,\cdot)$ follows the sign structure of $p$ over $\{1,\ldots,d+\alpha\}$. Moreover, since $B_{\alpha}^u$ and $B_{\alpha}^v$ realize $p_u$ and $p_v$, we have for every $k \in \{1, \ldots, d+\alpha\}$
\[
\begin{aligned}
\sum_{j=1}^d B_{\alpha}(j,k)&=\tau \sum_{j=1}^dB_{\alpha}^{u}(j,k)+(1-\tau)\sum_{j=1}^dB_{\alpha}^{v}(j,k)\\[0mm]
&=\tau p_{u}(k)+(1-\tau)p_{v}(k)\\[2mm]
&=\frac{\tau(k-u)+(1-\tau)(k-v)}{k-{k_{n_0}}}p(k)=p(k).
\end{aligned}
\]
We may therefore conclude that $B_{\alpha}$ realizes $p$ and the result holds.

\vspace{5pt}

\underline{\textit{Step III:}} Assume that the zeroes of $p$ are real. This means that $p$ has a decomposition
\[
p(x)=\frac{(-1)^{\alpha}}{\alpha!}\prod_{n=1}^{\alpha}(x-x_n),
\]
with $x_n\in\R$ not necessarily distinct.  We prove by induction on the \textit{non-integer degree} of $p$, which we define by
\[
\mathrm{ni}(p):=\#\{n:x_n\notin\Z\}.
\]
The case $\mathrm{ni}(p)=0$ was proved in \textit{Step II} above. So, assume that $\mathrm{ni}(p)\geq 1$. This means that there exists some $n_0\in\{1,...,\alpha\}$ so that $x_{n_0}\notin \Z$. Let
$u = \lfloor x_{n_0} \rfloor$ be the biggest integer smaller than $x_{n_0}$, and put $v =u+1$.
We now complete the proof exactly as in \textit{Step II}. That is, we define $p_u$ and $p_v$, as above, and note that their non-integer degree is $\mathrm{ni}(p)-1$. We then apply the induction assumption to find left inverses $B^u$ and $B^v$, respectively, that realise them. Finally, we take a convex combination of these to complete the proof.
\end{proof}

\subsection{Proof of the estimate (\ref{poly estimate})}

We are now ready to prove the estimate (\ref{poly estimate}). Recall that $\mathscr{P}_\alpha$ denotes the family of polynomials with leading coefficient equal to 1 and of degree equal to $\alpha$. Let $\widetilde{\mathscr{P}}_\alpha\subseteq\mathscr{P}_\alpha$ denote the subfamily of polynomials with real zeroes. It is straightforward to check that
\[
\inf_{p\in\mathscr{P}_\alpha} \left(\sum_{k=1}^{d+\alpha}|p(k)|^2\right)=\inf_{p\in\widetilde{\mathscr{P}}_\alpha} \left(\sum_{k=1}^{d+\alpha}|p(k)|^2\right).
\]
It is therefore enough to show that
\begin{equation}\label{poly estimate real}
\norm{ (a_j)_{j=1}^{d} }_{\ell^1} \leq
\frac{1}{\alpha!} \inf_{p\in\widetilde{\mathscr{P}}_\alpha} \left(\sum_{k=1}^{d+\alpha}|p(k)|^2\right)^\frac12
\norm{((\Delta_d^\alpha a)_j)_{j=1}^{d+\alpha}}_{\ell^2}.
\end{equation}
For this, observe that if $B_\alpha$ is a left inverse of $\Delta_d^\alpha$ then for $a\in\mathbb{C}^d$ we have
\[
a_j= \sum_{k=1}^{d+\alpha}B_\alpha(j,k)(\Delta_d^\alpha a)_k, \quad j=1,2,\dots,d.
\]
Therefore we get,
\begin{equation}\label{CS}
\begin{aligned}
\norm{ (a_j)_{j=1}^{d} }_{\ell^1} 
&=
\sum_{j=1}^{d} \Big| \sum_{k=1}^{d+\alpha}B_\alpha(j,k)(\Delta_d^\alpha a)_k \Big| \\[2mm]
&\leq \sum_{k=1}^{d+\alpha} |(\Delta_d^\alpha a)_k|  \Big( \sum_{j=1}^{d} |B_\alpha(j,k)| \Big) \\[2mm]
&\leq
\norm{((\Delta_d^\alpha a)_j)_{j=1}^{d+\alpha}}_{\ell^2} \left(\sum_{k=1}^{d+\alpha}\left(\sum_{j=1}^{d}|B_\alpha(j,k)|\right)^2\right)^\frac12,
 \end{aligned}
\end{equation}
where the Cauchy-Schwarz inequality was applied in the last step.
Let $p$ be a monic polynomial of degree $\alpha$ and whose zeroes are all real. By Theorem \ref{main thm} there exists a left inverse $B_{\alpha}$ which realizes $(-1)^{\alpha}p/\alpha !$ over $\{1,\ldots,d+\alpha\}$. Remark \ref{abs value equality} therefore implies that
\[
\sum_{j=1}^d |B_{\alpha}(j,k)|= \frac{1}{\alpha !}|p(k)|.
\]
We substitute this in (\ref{CS}) and take the infimum over all such monic polynomials to conclude (\ref{poly estimate real}).

\section{Proof of Theorem \ref{main_theorem} (part 2 of 2)}

In this section we complete the proof of Theorem \ref{main_theorem} in two steps. First, in Subsection \ref{subsection 4.1}, we survey the results which imply the identity
\begin{equation}\label{chebi estimate}
\min_{p\in\mathscr{P}_\alpha} \left(\sum_{k=1}^{d+\alpha}|p(k)|^2\right)^\frac12=\alpha!\sqrt{\binom{d+2\alpha}{2\alpha +1}\binom{2\alpha}{\alpha}^{-1}}.
\end{equation}
Combined with the estimate (\ref{poly estimate}), which was obtained at the end of the last section, this gives the estimate (\ref{the_main}) from  Theorem \ref{main_theorem}. Next, in Subsection \ref{subsection 4.2}, we show that the estimate (\ref{the_main}) is sharp.

\subsection{A proof for the identity (\ref{chebi estimate})}\label{subsection 4.1}

For $N\in\N$ denote
$
\mu_N:=\sum_{j=1}^N\delta_j,
$
where $\delta_j$ is the Dirac delta measure located at the integer $j$. With this notation,  the left hand side of (\ref{chebi estimate}) asks to minimize the norm $\|p\|_{L^2(\mu_{d+\alpha})}$ over all polynomials $p\in \mathscr{P}_{\alpha}$.

As is standard, a solution to this minimization problem may be found via the family of orthogonal polynomials in the space. The corresponding polynomials are the \textit{discrete Chebyshev polynomials} (see, e.g., \cite[p. 8]{Szego1975}). 

To introduce these functions, let $N\in\N$ and
consider the discrete Landau polynomials 
\begin{equation}\label{the L polynomial}
L^{(N)}_n(x):= r_n(x-1)r_n\big(x-(N+1)\big)\qquad n=0,...,N-1,
\end{equation}
where the $r_{n}$ are as defined in (\ref{r_a}).
With this, the discrete Chebyshev  polynomials are defined by
\begin{equation}\label{the T polynomial}
T^{(N)}_n :=(n!)\binom{2n}{n}^{-1}\Delta^nL^{(N)}_n,
\end{equation}
where $\Delta f(x)=f(x)-f(x-1)$. That is, $\Delta$ is an extension of the operator introduced in \eqref{discrete delta on Z} to functions defined over $\R$.

The following lemma records the main properties of the discrete Chebyshev  polynomials.

\begin{lemma}\label{chevilemma}
For $N\in\N$, the following hold.
\begin{itemize}
\item[(i)] For every $n\in \{0,...,N-1\}$, the polynomial $T^{(N)}_n$ is monic and of degree $n$.

\vspace{5pt}

\item[(ii)] The family $\{T^{(N)}_n\}_{n=0}^{N-1}$ is orthogonal in the space $L^2(\mu_N)$.

\vspace{5pt}

\item[(iii)] For every $n\in \{0,...,N-1\}$ we have
\[
\|T^{(N)}_n\|_{L^2(\mu_{N})}^2= (n!)^2\binom{2n}{n}^{-1}\binom{N+n}{2n+1}.
\]
\end{itemize}
\end{lemma}
\begin{proof}
Part (i) follows from applying $\Delta^{n}$ to a polynomial of degree $2n$ with leading coefficient $(n!)^{-2}$. Parts (ii) and (iii) are the formulas (2.8.2) and (2.8.3)
from \cite{Szego1975}, reformulated in the notations used here (note, in particular, that the operator $\Delta$ in \cite{Szego1975} is the forward difference operator $\Delta a_j = a_{j+1}-a_j$, whereas we use $\Delta$ to denote the backward difference operator).
\end{proof}

For brevity, put $\mu = \mu_{d+\alpha}$. An immediate consequence of Lemma \ref{chevilemma} is the observation that
\begin{equation}\label{minbychev}
\min_{p\in\mathscr{P}_\alpha}\|p\|_{L^2(\mu)}=\|T^{(d+\alpha)}_{\alpha}\|_{L^2(\mu)},
\end{equation}
To see why $T^{(d+\alpha)}_{\alpha}$ minimizes this norm, we first note that part (i) of Lemma \ref{chevilemma} implies that the family $\{T^{(d+\alpha)}_{n}\}_{n=0}^{\alpha-1}$ spans the space of all polynomials of degree less than $\alpha$. Part (ii) of Lemma \ref{chevilemma} now implies that $T^{(d+\alpha)}_{\alpha}$ is orthogonal to this space in $L^2(\mu)$.
Next, let $p\in\mathscr{P}_\alpha$. Since both $T^{(d+\alpha)}_{\alpha}$ and $p$ are monic, the polynomial $q=p-T^{(d+\alpha)}_{\alpha}$ is of degree at most $\alpha-1$. By the above remark we may conclude that  $T^{(d+\alpha)}_{\alpha}$ is orthogonal to it. This implies that,
\[
\|p\|^2_{L^2(\mu)}=\|T^{(d+\alpha)}_{\alpha}\|^2_{L^2(\mu)}+\|q\|^2_{L^2(\mu)}\geq\|T^{(d+\alpha)}_{\alpha}\|^2_{L^2(\mu)}.
\]
With this (\ref{minbychev}) is established.
By substituting part (iii) of Lemma \ref{chevilemma} into (\ref{minbychev}), we obtain the identity (\ref{chebi estimate}).

\subsection{Sharpness}\label{subsection 4.2}

Our goal in this subsection is to show that the constant in the inequality (\ref{the_main}) is sharp. That is, we find $(a_j)_{j=1}^d$ so that
	\begin{equation}\label{sharpness}
		\norm{ (a_j)_{j=1}^{d} }_{\ell^1} = \sqrt{\binom{d+2\alpha}{2\alpha +1}\binom{2\alpha}{\alpha}^{-1}} \norm{((\Delta_d^\alpha a)_j)_{j=1}^{d+\alpha}}_{\ell^2}.
	\end{equation}
For this, recall the discrete Landau polynomials  $L^{(N)}_n$ defined in (\ref{the L polynomial}) and, for $j \in \{1,\ldots, d\}$, put
\begin{align*}
	a_j &:=L^{(d+\alpha)}_{\alpha}(j).
\end{align*}
To establish that the identity (\ref{sharpness}) holds,  we first observe that the polynomial $L^{(d+\alpha)}_{\alpha}$ vanishes over the sets $\{-\alpha+1,...,0\}$ and $\{d+1,...,d+\alpha\}$. Consequently, we have for every $j \in \{1,\ldots, d\}$ that
\[
(\Delta^{\alpha}_d a)_j=\Delta_d^{\alpha}L^{(d+\alpha)}_{\alpha}(j)=\frac{1}{\alpha!}\binom{2\alpha}{\alpha}T^{d+\alpha}_{\alpha}(j).
\]
Part (iii) of Lemma \ref{chevilemma} now implies that
\[
\norm{\Delta_d^\alpha a}_{\ell^2}=\sqrt{\binom{d+2\alpha}{2\alpha +1}\binom{2\alpha}{\alpha}}.
\]
Therefore, to establish the identity (\ref{sharpness}), it remains to prove the following proposition.

\begin{proposition} With the notations above, we have
\begin{equation}\label{ell_1 norm of gamma}
\norm{ a }_{\ell^1}=\binom{d+2\alpha}{2\alpha +1}.
\end{equation}
\end{proposition}
\begin{proof}
We first note that, due to the sign structure of $r_{\alpha}$, we have
\[
\norm{ a }_{\ell^1}=(-1)^{\alpha}\sum_{j=1}^dr_{\alpha}(j-1)r_{\alpha}\big(j-(d+\alpha+1)\big).
\]

Next, we observe that by \eqref{Delta r_alpha}, we have $r_\alpha = \Delta^\alpha r_{2\alpha}$. Applying this above, we obtain
\[
\norm{ a }_{\ell^1}=(-1)^{\alpha}\sum_{j=1}^d\Delta^{\alpha}r_{2\alpha}(j-1)r_{\alpha}\big(j-(d+\alpha+1)\big).
\]
Applying summation by parts several times,  keeping in mind that the boundary terms vanish, we get
\[
\norm{ a }_{\ell^1}=\sum_{j=1}^dr_{2\alpha}(j-1)\Delta^{\alpha} r_{\alpha}(j-(d+1)).
\]
Finally,   $\Delta^{\alpha} r_{\alpha}=r_{0}=1$, and so we are left with
\[
\begin{aligned}
\norm{ a }_{\ell^1}&=\sum_{j=1}^dr_{2\alpha}(j-1)
=\sum_{j=1}^d\binom{2\alpha+j-1}{2\alpha}\\
&=\sum_{j=2\alpha}^{2\alpha+d-1}\binom{j}{2\alpha}
=\binom{d+2\alpha}{2\alpha+1},
\end{aligned}
\]
where the last equality is due to the \textit{hockey-stick} binomial identity. This completes the proof.
\end{proof}

\bibliographystyle{plain}

\bigskip

\noindent
\textsc{Centre for Mathematical Sciences, Lund University,
Box 118, 221 00 Lund, Sweden}

\noindent
E-mail address: \texttt{raul.hindov@math.lu.se}

\medskip

\noindent
\textsc{Georgia Institute of Technology, Atlanta, USA}

\noindent
E-mail address: \texttt{shahaf.nitzan@math.gatech.edu}

\medskip

\noindent
\textsc{Centre for Mathematical Sciences, Lund University,
Box 118, 221 00 Lund, Sweden}

\noindent
E-mail address: \texttt{jan-fredrik.olsen@math.lu.se}

\medskip

\noindent
\textsc{Centre for Mathematical Sciences, Lund University,
Box 118, 221 00 Lund, Sweden}

\noindent
E-mail address: \texttt{eskil.rydhe@math.lu.se}

\end{document}